\documentclass[11pt,article]{amsart}
\usepackage{amssymb}
\usepackage{amsfonts}
\usepackage{amsbsy}
\usepackage{latexsym}

\usepackage{graphicx}

\usepackage{amsmath}

\newtheorem{theorem}{Theorem}[section]

\newtheorem{proposition}[theorem]{Proposition}
\newtheorem{corollary}[theorem]{Corollary}

\theoremstyle{definition}
\newtheorem{definition}[theorem]{Definition}
\newtheorem{example}[theorem]{Example}

\theoremstyle{remark}
\newtheorem{remark}[theorem]{Remark}

\numberwithin{equation}{section}

\newcommand\N{\mathbb{N}}

\newcommand\C{\mathbb{C}}
\newcommand\J{\mathbb{J}}

\newcommand\cH{\mathcal{H}}
\newcommand\cF{\mathcal{F}}
\newcommand\cK{\mathcal{K}}
\newcommand\cG{\mathcal{G}}
\newcommand\cB{\mathcal{B}}
\newcommand\cD{\mathcal{D}}
\newcommand\cM{\mathcal{M}}

\def\sideremark#1{\ifvmode\leavevmode\fi\vadjust{\vbox
to0pt{\vss \hbox to 0pt{\hskip\hsize\hskip1em
\vbox{\hsize2cm\tiny\raggedright\pretolerance10000
\noindent#1\hfill}\hss}\vbox to8pt{\vfil}\vss}}}

\newcommand{\ip}[2]{\left\langle\,#1 ,  #2 \, \right\rangle}

\begin{document}


\title{Dilations of frames,  operator valued measures and bounded linear maps}

\author{Deguang Han}
\address{Department of Mathematics, University of Central
Florida, Orlando, USA} \email{deguang.han@ucf.edu}

\author{David R. Larson}
\address{Department of Mathematics, Texas A\&M University, College Station, USA}
\email{larson@math.tamu.edu}

\author{Bei Liu}
\address{Department of Mathematics, Tianjin University of Technology, Tianjin
300384, P.R. China} \email{beiliu1101@gmail.com}

\author{Rui Liu}
\address{Department of Mathematics and LPMC, Nankai University, Tianjin 300071, P.R. China}
\address{Department of Mathematics, Texas A\&M University, College Station, USA}
\email{ruiliu@nankai.edu.cn}

\thanks{{\it Acknowledgements}: The authors were all participants in the NSF funded Workshop
in Analysis and Probability at Texas A\&M University. The first
author acknowledges partial support by a grant from the NSF. The
third and fourth authors received partial support from the NSFC}

\date{}


\begin{abstract}

We will give an outline of the main results in our recent AMS
Memoir, and include some new results, exposition and open problems.
In that memoir we developed a general dilation theory for operator
valued measures acting on Banach spaces where operator-valued
measures (or maps) are not necessarily completely bounded.  The main
results state that  any operator-valued measure, not necessarily
completely bounded, always has a dilation to a projection-valued
measure acting on a Banach space, and every bounded linear map,
again not necessarily completely bounded, on a Banach algebra has a
bounded homomorphism dilation acting on a Banach space.   Here the
dilation space often needs to be a Banach space even if the
underlying space is a Hilbert space, and the projections are
idempotents that are not necessarily self-adjoint. These results
lead to some new connections between frame theory and operator
algebras, and some of them can be considered as part of the
investigation  about ``noncommutative" frame theory.
\end{abstract}

\maketitle


\section{Introduction}

Frame theory belongs to the area of applied harmonic analysis, but
its underpinnings involve large areas of functional analysis
including operator theory, and deep connections with the theory of
operator algebras on Hilbert space.  For instance, recently the
Kadison-Singer ``Extension of pure states on von Neumann algebras"
problem has been solved and its solution is known to have wide
ramifications in frame theory due mainly to the research and
excellent exposition of Casazza and of Weaver.  The purpose of the
present article is to give a good exposition of some recent work of
the authors that establishes some deep connections between frame
theory on the one hand and operator-valued measures and maps between
von Neumann algebras on the other hand. While our work has little or
nothing to do directly with the above-mentioned extension of pure
states problem, it represents a separate instance of a rather deep
connection between frame theory and operator algebras, and this is
the point of this article.  We show that it may have something to do
with another problem of Kadison: the ``similarity problem".  At the
least it indicates the possibility of another possible approach to
that problem.  And it does show that the ideas implicit in frame
theory belong to the underpinnings of a significant part of modern
mathematics.

Before discussing our results from \cite{HLLL1} we need an
exposition of some of the preliminaries leading up to them.

\bigskip

\section{Frames, Framings, and Operator Valued Measures}

A \textit{frame} $\mathcal{F}$ for a Hilbert space $\cH$ is a
sequence of vectors $\{x_n\} \subset \cH$ indexed by a countable
index set $\J$ for which there exist constants $0<A\leq B<\infty$
such that, for every $x \in \cH$,

\begin{equation}
\label{frame} A\|x\|^2 \leq \sum_{n\in\J} |\ip{x}{x_n}|^2 \leq
B\|x\|^2
\end{equation}
The optimal constants are known as the upper and lower \textit{frame
bounds}. A
 frame is called  \textit{tight} if $A=B$,  and
is called a \textit{Parseval frame} if $A=B=1$. If we only require
that a sequence $\{x_{n}\}$ satisfies the upper bound condition in
(\ref{frame}), then $\{x_{n}\}$ is called a {\it Bessel sequence}. A
frame which is a basis is called a Riesz basis. Orthonormal bases
are special cases of Parseval frames. A Parseval frame $\{x_n\}$ for
a Hilbert space $\cH$ is an orthonormal basis if and only if each
$x_n$ is a unit vector.

For a Bessel sequence $\{x_{n}\}$, its {\it analysis operator}
$\Theta$ is a bounded linear operator from $\cH$ to
$\ell^{2}(\mathbb{N})$ defined by
\begin{equation}\label{frameDef}
\Theta x = \sum_{n\in \N}\ip{x}{x_{n}}e_{n},
\end{equation}
 where $\{e_{n}\}$ is the standard
orthonormal basis for $\ell^{2}(\mathbb{N})$. It is easily verified
that
$$
\Theta^{*}e_{n} = x_{n}, \ \ \forall n\in\N
$$
The Hilbert space adjoint $\Theta^{*}$ is called the {\it synthesis
operator} for $\{x_{n}\}$. The positive operator
$S:=\Theta^{*}\Theta:\cH \rightarrow \cH$ is called the {\it frame
operator}, or sometimes the {\it Bessel operator} if the Bessel
sequence is not a frame, and we have
\begin{equation}\label{frameOp}
Sx = \sum_{n\in \N}\ip{x}{ x_{n}}x_{n},  \ \ \ \forall x\in \cH.
\end{equation}

A sequence $\{x_{n}\}$ is a  frame for $\cH$ if and only if its
analysis operator $\Theta$ is bounded, injective and has closed
range, which is, in turn, equivalent to the condition that the frame
operator $S$ is bounded and invertible. In particular, $\{x_{n}\}$
is a  Parseval frame for $\cH$ if and only if $\Theta$ is an
isometry or equivalently  if $S = I$.

Let $S$ be the frame operator for a frame $\{x_{n}\}$. Then the
lower frame bound is $1/||S^{-1}||$ and the upper frame bound is
$||S||$. From (\ref{frameOp}) we obtain the {\it reconstruction
formula (or frame decomposition)}: $$x = \sum_{n\in \N}\ip{x}{
S^{-1}x_{n}}x_{n}, \ \ \forall x\in \cH$$ or equivalently $$ x =
\sum_{n\in \N}\ip{x}{x_{n}}S^{-1}x_{n}, \ \ \forall x\in \cH. $$

The frame $\{S^{-1}x_{n}\}$ is called the {\it canonical or standard
dual} of $\{x_{n}\}$. In the case that $\{x_{n}\}$ is a Parseval
 frame for $\cH$, we have that $S = I$ and hence
$ x = \sum_{n\in \N}\ip{x}{x_{n}}x_{n}, \ \ \forall x\in \cH. $ More
generally, if a Bessel sequence $\{y_{n}\}$ satisfies
$x = \sum_{n\in \N}\ip{x}{y_{n}}x_{n}, \ \ \forall x\in \cH, $
where the convergence is in norm of $\cH$, then $\{y_{n}\}$ is
called an {\it alternate dual } of $\{x_{n}\}$. (Then $\{y_{n}\}$ is
also necessarily a frame.) The canonical and alternate duals are
often simply referred to as {\it duals}, and $\{x_{n}, y_{n}\}$ is
called a {\it dual frame pair.} It is a well-known fact that that a
frame $\{x_{n}\}$ is a Riesz basis if and only if $\{x_{n}\}$  has a
unique dual frame.

There is a geometric interpretation of Parseval frames and general
frames. Let $P$ be an orthogonal projection from a Hilbert space
$\cK$ onto a closed subspace $\cH$, and let $\{u_{n}\}$ be a
sequence in $\cK$. Then $\{Pu_{n}\}$ is called the {\it{orthogonal
compression}} of $\{u_{n}\}$ under $P$, and correspondingly
$\{u_{n}\}$ is called an {\it{orthogonal dilation}} of $\{Pu_{n}\}$.
We first observe that if $\{u_n\}$ is a  frame for $\cK$, then
$\{Pu_n\}$ is a  frame for $\cH$ with frame bounds at least as good
as those of $\{u_n\}$ (in the sense that the lower frame cannot
decrease and the upper bound cannot increase). In particular,
$\{Pu_n\}$ is a Parseval frame for $\cH$ when $\{u_{n}\}$ is an
orthonormal basis for $\cK$, i.e., every orthogonal compression of
an orthonormal basis (resp. Riesz basis) is a Parseval frame (resp.
frame) for the projection subspace. The converse is also true: every
frame can be orthogonally dilated to a Riesz basis, and every
Parseval frame can be dilated to an orthonormal basis. This was
apparently first shown explicitly by Han and Larson in Chapter 1 of
\cite{HL}. There, with appropriate definitions it had an elementary
two-line proof. And as noted by several authors, it can be
alternately derived by applying the Naimark (Neumark) Dilation
theorem for operator valued measures by first passing from a frame
sequence to a natural discrete positive operator-valued measure on
the power set of the index set. So it is sometimes referred to as
the Naimark dilation theorem for frames.  In fact, this is the
observation that inspired much of the work in \cite{HLLL1}.

For completeness we formally state this result:
\begin{proposition} \cite{HL} \label{prop:orthog} Let $\{x_{n}\}$ be a
sequence in a Hilbert space $\cH$. Then
\begin{enumerate}

\item  $\{x_n\}$ is a Parseval frame for $\cH$ if and only if
there exists a Hilbert space $\cK \supseteq \cH$ and an orthonormal
basis $\{u_n\}$ for $\cK$ such that $x_n = Pu_n$, where $P$ is the
orthogonal projection from $\cK$ onto $\cH$.

\item $\{x_n\}$ is a  frame for  $\cH$ if and only if there exists
a Hilbert space $\cK \supseteq \cH$ and a Riesz basis $\{v_n\}$ for
$\cK$ such that $x_n = Pv_n$, where $P$ again is the orthogonal
projection from $\cK$ onto $\cH$.
\end{enumerate}
\end{proposition}

The above dilation result was later generalized in \cite{CHL} to
dual frame pairs.

\begin{theorem} \label{ch2-dualpair-dilation}
Suppose that  $\{x_n\}$ and $\{y_{n}\}$ are two frames for a Hilbert
space $\cH$. Then the following are equivalent:
\begin{enumerate}
\item $\{y_{n}\}$ is a dual for $\{x_{n}\}$;

\item  There exists a Hilbert space $\cK \supseteq \cH$ and a Riesz
basis $\{u_n\}$ for $\cK$ such that $x_n = Pu_n$, and $y_{n} =
Pu_{n}^{*}$, where $\{u_{n}^{*}\}$ is the (unique) dual of the Riesz
basis $\{u_{n}\}$ and $P$  is the orthogonal projection from $\cK$
onto $\cH$.
\end{enumerate}
\end{theorem}

As in  \cite{CHL}, a {\it framing} for a Banach space $X$ is a pair
of sequences $\{x_i,y_i\}$ with $\{x_i\}$ in $X$, $\{y_i\}$ in the
dual space $X^*$ of $X$, satisfying the condition that
\[x=\sum_i\langle x, y_i\rangle x_i,\]
where this series converges unconditionally for all $x\in X.$

The definition of a framing is a natural generalization of the
definition of a dual frame pair. Assume that $\{x_i\}$ is a  frame
for $\cH$ and $\{y_i\}$ is a dual frame for $\{x_i\}$. Then
$\{x_i,y_i\}$ is clearly a framing for $\cH.$ Moreover, if
$\alpha_i$ is a sequence of non-zero constants, then $\{\alpha_ix_i,
\bar{\alpha}^{-1}_iy_i\}$ (called a rescaling of the pair) is also a
framing, although it is easy to show
that it need not be a pair of frames, even if $\{\alpha_ix_i\}$,
$\{\bar{\alpha}^{-1}_iy_i\}$ are bounded sequence.

We recall that a sequence $\{z_i\}$ in a Banach space $Z$ is called a \emph{Schauder basis} (or
just a \emph{basis}) for $Z$ if for each $z\in Z$ there is a unique sequence of scalars $\{\alpha_i\}$
so that $z = \sum_i \alpha_i z_i$. The unique elements $z_i^*\in Z^*$ satisfying
\begin{equation}\label{eq:2.2} z=\sum_i z_i^*(z)z_i, \end{equation}
for all $z\in Z,$ are called the \emph{dual} (or \emph{biorthogonal}) \emph{functionals} for $\{z_i\}$.
If the series in (\ref{eq:2.2})
converges unconditionally for every $z\in Z$, we call $\{z_i, z_i^*\}$ an \emph{unconditional
basis} for $Z$.

We also have an \emph{unconditional basis constant} for an unconditional
basis given by:
\[ \mathrm{UBC}(z_i) = \sup\{ \|\sum_i b_i z_i\| : \| \sum_i a_i z_i \| =1, |b_i|\le |a_i|, \forall i \}. \]
If $\{z_i, z_i^*\}$ is an unconditional basis for $Z$, we can define an equivalent norm on $Z$ by:
\[ \|\sum_i a_i z_i\|_1 = \sup\{ \|\sum_i b_ia_iz_i\| : |b_i|\le 1, \forall i \}. \]
Then $\{z_i, z_i^*\}$ is an unconditional basis for $Z$ with $\mathrm{UBC}(z_i) = 1$. In this case we
just call $\{z_i\}$ a $1$-unconditional basis for $Z$.

\begin{definition}\label{de:31} \cite{CHL}
A sequence $\{x_i\}_{i\in \N}$ in a Banach space $X$ is a {\it
projective frame for $X$} if there is a Banach space $Z$ with an
unconditional basis $\{z_i, z^*_i\}$ with $X\subset Z$ and a (onto)
projection $P: Z\to X$ so that $Pz_i = x_i$ for all $i \in\N$. If
$\{z_i\}$ is a 1-unconditional basis for $Z$ and $\|P\| = 1$, we
will call $\{x_i\}$ a {\it projective Parseval} frame for $X$.
\end{definition}
In this case, we have for all $x \in X$ that
\[x =\sum_i\langle
x,z^*_i\rangle z_i = Px =\sum_i\langle x,z^*_i\rangle Pz_i
=\sum_i\langle x,z^*_i\rangle x_i,\] and this series converges
unconditionally in X. So this definition recaptures the
unconditional convergence from the Hilbert space definition.

We note that there exist projective frames in the sense of
Definition \ref{de:31} for an infinite dimensional Hilbert space
that fail to be  frames.  We think they occur in abundance, but
specific examples are hard to prove. A concrete example is contained
in \cite[Chapter 5]{HLLL1}.

\begin{definition} \label{de:32} \cite{CHL}
A framing model is a Banach space $Z$ with a fixed unconditional
basis $\{e_i\}$ for $Z.$ A framing modeled on $(Z, \{e_i\}_{i\in\N}
)$ for a Banach space $X$ is a pair of sequences $\{y_i\}$ in $X^*$
and $\{x_i\}$ in $X$ so that the operator $\theta: X\to Z$ defined
by
\[\theta u =\sum_{i\in\N}\langle u,y_i\rangle e_i,\]
 is an into isomorphism and $\Gamma: Z \to X$ given by
 \[\Gamma(\sum_{i\in\N}
a_ie_i)=\sum_{i\in\N}a_ix_i\] is bounded and $\Gamma \theta=I_X$.
\end{definition}

In this setting, $\Gamma$ becomes the reconstruction operator for
the frame.  The following result due to Casazza, Han and Larson
\cite{CHL}
 shows that these three methods for defining a frame on a
Banach space are really the same.
\begin{proposition}\label{pr:33}
Let $X$ be a Banach space and $\{x_i\}$ be a sequence of elements of
$X.$ The following are equivalent:
\begin{enumerate}
\item[(1)]$\{x_i\}$ is a projective frame for $X$.

\item[(2)]There exists a sequence $y_i\in X^*$  so that $\{x_i,
y_i\}$ is a framing for $X.$

\item[(3)]There exists a sequence $y_i\in X^*$  and a framing
model $(Z, \{e_i\})$ so that $\{x_i, y_i\}$ is a framing modeled on
$(Z, \{e_i\}).$
\end{enumerate}
\end{proposition}

The proof of the implication from $(1)$ to $(2)$ is trivial: If $\{z_{i}\}$ is an unconditional basis for a Banach space $Z$ and $P$ is a bounded projection from $Z$ to a closed subspace $X$ with $x_{i} = Pe_{i}$, then $(x_{i}, y_{i})$ is a framing for $X$, where $y_{i} = P^{*}z_{i}^{*}$ and $\{z_{i}^{*}\}$ is the (unique) dual basis of $\{z_{i}\}$. One of the main contributions of  paper \cite{CHL} was to show that every framing can be obtained in this way.

\begin{theorem} [Corollary 4.7 of \cite{CHL}]  Suppose that $\{x_i, y_i\}$ is a framing for
$X$. Then there exist a Banach space $Z$  containing $X$ and an unconditional basis
$\{z_i, z^*_i\}$  for $Z$  such that $x_i = Pz_{i}$ and $y_{i} = P^{*}z_{i}^{*}$, where $P$ is a bounded projection from $Z$ onto $X$.
\end{theorem}

The definition of (discrete) frames has a natural generalization.

\begin{definition}
Let $\cH$ be a separable Hilbert space and $\Omega$ be a
$\sigma$-locally compact ($\sigma$-compact and locally compact)
Hausdorff space endowed with a positive Radon measure $\mu$ with
$\mbox{supp}(\mu) =\Omega$. A weakly continuous function
$\mathcal{F}:\Omega\to \cH$ is called a \emph{continuous frame} if
there exist constants $0 < C_1\le C_2 < \infty$ such that

\begin{equation*}
 C_1\|x\|^2 \le \int_\Omega |\langle x,\mathcal{F}(\omega)
\rangle|^2 d \,\mu(\omega) \le C_2\|x \|^2, \quad \forall\, x\in
\cH.
\end{equation*}

\end{definition}

If $C_1 = C_2$ then the frame is called \emph{tight}. Associated to
$\mathcal{F} $ is the frame operator $S_\mathcal{F} : \cH \to \cH$
defined in the weak sense by
\begin{equation*} \quad \langle S_\cF  (x),y\rangle
:= \int_\Omega \langle x, \mathcal{F} (\omega)\rangle\cdot\langle
\cF (\omega),y\rangle d \,\mu(\omega).
\end{equation*}

It follows from the definition that $S_\cF $ is a bounded, positive,
and invertible operator. We define the following transform
associated to $\cF $,
\begin{eqnarray*}
V_\cF  : \cH\to L^2(\Omega, \mu), \quad  V_\cF  (x)(\omega) :=
\langle x, \cF (\omega)\rangle.
\end{eqnarray*}
This operator is called the {\it analysis operator} in the
literature and its adjoint operator is given by
\begin{eqnarray*}
V_\cF ^* : L^2(\Omega, \mu) \to \cH, \quad \langle V_\cF
^*(f),x\rangle :=\int_{\Omega} f(\omega)\langle\cF (\omega),x\rangle
d\,\mu(\omega).
\end{eqnarray*}
Then we have $ S_\cF = V_\cF ^{*} V_\cF$, and
\begin{eqnarray} \label{cfdual}
\langle x,y\rangle = \int_\Omega \langle x, \cF
(\omega)\rangle\cdot\langle \cG (\omega),y\rangle d \,\mu(\omega),
\end{eqnarray}
where $\cG (\omega) := S_\cF ^{-1}\cF (\omega)$ is the {\it standard
dual} of $\cF$. A weakly continuous function $\cF :\Omega\to \cH$ is
called \emph{Bessel} if there exists a positive constant $C$ such
that
\begin{eqnarray*}
\int_\Omega |\langle x,\cF (\omega) \rangle|^2 d \,\mu(\omega) \le
C\|x \|^2, \quad \forall\, x\in \cH.
\end{eqnarray*}

It can be easily shown  that if $\cF :\Omega\to \cH$ is  Bessel,
then it is a frame for $\cH$ if and only if there exists a Bessel
mapping $\cG$ such that the reconstruction formula (\ref{cfdual})
holds. This $\cG$ may not be the standard dual of $\cF$. We will
call $(\cF, \cG)$ a {\it dual pair}.

A discrete frame is a Riesz basis if and only if its analysis
operator is surjective. But for a continuous frame $\cF$, in general
we don't have the dilation space to be $L^{2}(\Omega, \mu)$. In fact, this
could happen only when $\mu$ is  purely atomic. Therefore there is
no Riesz basis type dilation theory for continuous frames (however,
we will see later that in contrast the induced operator-valued
measure does have a projection valued measure dilation). The
following modified dilation theorem is due to Gabardo and Han
\cite{GH}:

\begin{theorem} Let $\cF$ be a $(\Omega, \mu)$-frame for
$\cH$ and $\cG$ be one of its duals. Suppose that both $V_\cF(\cH)$
and $V_\cG (\cH)$ are contained in the range space $\cM$ of the
analysis operator for some $(\Omega, \mu)$-frame. Then there is a
Hilbert space $\cK\supset \cH$ and a $(\Omega, \mu)$-frame
$\tilde{\cF}$ for $\cK$ with  $P\tilde{\cF} = \cF$, $P\tilde{\cG} =
\cG$ and $V_{\tilde{\cF}}(\cH) = \cM$, where $\tilde{\cG}$ is the
standard dual of $\tilde{\cF}$ and $P$ is the orthogonal projection
from $\cK$ onto $\cH$.
\end{theorem}

Let $\Omega$ be a compact Hausdorff space, and let $\cB$ be the
$\sigma$-algebra of all the Borel subsets of $\Omega$. A
$B(\cH)$-valued measure on $\Omega$ is a mapping $E:\cB\to B(\cH)$
that is weakly countably additive, i.e., if $\{B_i\}$ is a countable
collection of disjoint Borel sets with union $B$, then
\[\langle E(B)x,y\rangle=\sum_i\langle E(B_i)x,y\rangle\]
holds for all $x,y$ in $\cH$. The measure is called {\it bounded}
provided that
\[\sup\{\|E(B)\|:B\in \cB\}<\infty,\]
and we let $\|E\|$ denote this supremum. The measure is called
{\it regular} if for all $x,y$ in $\cH$, the complex measure given
by
\begin{eqnarray}\label{eq:w141}
 \mu_{x,y}(B)=\langle E(B)x,y \rangle
\end{eqnarray} is regular.

Given a regular bounded $B(\cH)$-valued measure $E$, one obtains a
bounded, linear map
\[\phi_E:C(\Omega)\to B(\cH)\] by
\begin{eqnarray}\label{eq:w142}
\langle \phi_E(f)x,y \rangle=\int_\Omega f\, d\, \mu_{x,y}.
\end{eqnarray}

Conversely, given a bounded, linear map $\phi:C(\Omega)\to B(\cH)$,
if one defines regular Borel measures $\{\mu_{x,y}\}$ for each $x,y$
in $\cH$ by the above formula (\ref{eq:w142}), then for each Borel
set $B$, there exists a unique, bounded operator $E(B)$, defined by
formula (\ref{eq:w141}), and the map $B\to E(B)$ defines a bounded,
regular $B(\cH)$-valued measure. There is a one-to-one
correspondence between the bounded, linear maps of $C(\Omega)$ into
$B(\cH)$ and the regular bounded $B(\cH)$-valued measures. Such
measures are called
\begin{enumerate}
\item[(i)] {\it spectral}  if $E(B_1\cap B_2)=E(B_1)\cdot E(B_2),$ \item[(ii)] {\it positive} if $E(B)\geq 0$,
\item[(iii)] {\it self-adjoint} if $E(B)^*=E(B),$
\end{enumerate}
for all Borel sets $B,B_1$ and $B_2$.  Note that if $E$ is spectral
and self-adjoint, then $E(B)$ must be an orthogonal projection for
all $B\in \mathcal{B}$, and hence $E$ is positive.


In the commutative $C^*$ theory, compactness is usually used as
above because when viewing a unital $C^*$-algebra as $C(\Omega)$
there is no loss in generality in taking $\Omega$ to be compact,
because if needed it can be taken to be $\beta\Omega$ -- the
Stone-Cech compactification of $\Omega$.  This is because the
$C^*$-algebras $C(\Omega)$ and $C(\beta\Omega)$ are $*$-isomorphic.
Having $\Omega$ compact makes the integration theory representation
of linear maps and the connection between linear maps on $C(\Omega)$
and operator valued measures very  elegant.

But in our theory, the basic connection to frame theory is
essentially lost if we replace the index set of the frame with its
Stone-cech compactification. In the continous frame case it is more
natural to assume $\Omega$ is $\sigma$-locally compact (as in
Definition 2.7), and in the general dilation theory we need to use
the {\it general} measurable space setting (as in Definition 3.2) to
preserve our basic connections with the frame theory.


Both discrete and continuous framings induce operator valued measures in a natural way.

\begin{example} Let $\{x_i\}_{i\in\J}$ be a frame for a separable Hilbert
space $\cH.$ Let $\Sigma$ be the $\sigma$-algebra of all subsets of
$\J$. Define the mapping
$$E: \Sigma\to B(\cH), \quad
E(B)=\sum_{i\in B} x_i\otimes x_i$$ where $x\otimes y$ is the
mapping on $\cH$ defined by $(x\otimes y)(u)=\langle u,y\rangle x.$
Then $E$ is a regular, positive $B(\cH)$-valued measure.

Similarly, suppose that $\{x_i,y_i\}_{i\in\J}$ is a non-zero framing
for a separable Hilbert space $\cH.$ Define the mapping
$$E:\Sigma\to B(\cH), \quad
E(B)=\sum_{i\in B} x_i\otimes y_i, $$ for all $B\in \Sigma$. Then
$E$ is a $B(\cH)$-valued measure.
\end{example}

\begin{example}
Let $X$ be a Banach space and $\Omega$ be a $\sigma$-locally compact
Hausdorff space. Let $\mu$ be a Borel measure on $\Omega$. A
continuous framing on $X$ is a pair of maps $(\mathcal{F},\mathcal{G}),$
\[\mathcal{F}:\Omega\to X, \quad\mathcal{G}:\Omega\to X^*,\]
such that the equation
\begin{eqnarray*}
\left\langle E_{(\mathcal{F},\mathcal{G})}(B)x,y\right\rangle
=\int_B\langle x,\mathcal{G}(\omega)\rangle\langle
\mathcal{F}(\omega),y\rangle\,d\mu(\omega)
\end{eqnarray*} for $x\in X$, $y\in X^*$, and $B$ a Borel subset of $\Omega$,
defines an operator-valued probability measure on $\Omega$ taking
value in $B(X)$.  In particular, we
require the integral on the right to converge for each
$B\subset\Omega.$ We have
\begin{eqnarray}\label{eq:c71}
E_{(\mathcal{F},\mathcal{G})}(B)=\int_B
\mathcal{F}(\omega)\otimes\mathcal{G}(\omega)\,d E(\omega)
\end{eqnarray}
where the integral converges in the sense of Bochner. In particular,
since $E_{(\mathcal{F},\mathcal{G})}(\Omega)=I_X,$ we have for any
$x\in X$ that
\[\langle x, y\rangle =\int_\Omega\langle x,\mathcal{G}(\omega)\rangle\langle\mathcal{F}(\omega), y\rangle\,d E(\omega).\]
\end{example}

%
\begin{remark} We point out that there exists an operator space with
a (finite dimensional) projection valued (purely atomic)  probability measure that does not admit a framing.
Let $X$ be the space of all compact operators $T$ on $\ell_2$ which have a triangular representing matrix
with respect to the unit vector basis, i.e. \[T e_n=\sum_{m=1}^n a_{n,m} e_m\] for all $n\in\N$.
Let $X_n$ be the subspace of $X$ consisting of those $T\in X$ such that $T e_j=0$ for $j\neq n$
(i.e. for which $a_{j, m}=0$ unless $j=n$). It is clear that $X_n$ is isometric to $\ell_2^n$, $n=1, 2, \ldots$.
Moreover, it is trivial to check that $\{X_n\}_{n=1}^\infty$ forms an unconditional finite dimensional decomposition of $X$
which naturally induces a projection valued probability measure.
Let $P_n$ be the canonical projection from $X$ onto $X_n$ satisfy:
\begin{enumerate}
\item[(i)] $\dim(P_n(X))=\dim(X_n)=n$ for all $n\in\N$;
\item[(ii)] $P_nP_m=P_mP_n=0$ for any $n\neq m\in\N$;
\item[(iii)] $x=\sum_{n=1}^\infty P_n(x)$ for every $x\in X$.
\end{enumerate}
Let $\Sigma$ be the $\sigma$-algebra of all subsets of $\N$. Define
$E:\Sigma\to B(X)$ by $E(B)=\sum_{n\in B} P_n.$
Then $E$ is a projection valued probability measure with $\dim(E(\{n\}))=n$.
Nevertheless, it follows from the results of \cite{GL} that
$X$ does not have an unconditional basis and it is not even complemented in a space with an unconditional basis. Thus, by Proposition \ref{pr:33},
$X$ does not have a framing.
\end{remark}

Let $\mathcal {A}$ be a unital $C^*$- algebra. An
operator-valued linear map $\phi:\mathcal {A}\to B(\cH)$ is said to
be {\it positive} if $\phi(a^{*}a) \geq 0$ for every $a\in
\mathcal{A}$, and it is called {\it completely positive} (cp for abbreviation) if for
every $n$-tuple
 $a_{1}, ... , a_{n}$ of elements in $\mathcal{A}$, the matrix $(\phi(a_{i}^{*}a_{j}))$ is positive in the usual sense
 that for every $n$-tuple of vectors $\xi_{1}, ... ,
\xi_{n} \in \cH$, we have
\begin{eqnarray}
\sum_{i, j=1}^{n}\langle \phi(a_{i}a_{j}^{*})\xi_{j}, \xi_{i}\rangle
\geq 0
\end{eqnarray}
or equivalently, $(\phi(a_i^*a_j))$ is a positive operator on the
Hilbert space $\cH\otimes\C^n$ (\cite{Pa}).

Let $\mathcal {A}$ be a $C^*$- algebra. We use  $M_n$ to denote the
set of all $n\times n$ complex matrices, and  $M_n(\mathcal {A})$ to
denote the set of all $n\times n$ matrices with entries from
$\mathcal {A}.$
Given two $C^*$-algebras $\mathcal {A}$ and $\mathcal {B}$ and a map
$\phi:\mathcal {A}\to\mathcal {B}$, obtain maps $\phi_n:M_n(\mathcal
{A})\to M_n(\mathcal {B})$ via the formula
\[\phi_n((a_{i,j}))=(\phi(a_{i,j})).\]
The map $\phi$ is called \textit{completely bounded} (cb for abbreviation) if $\phi$ is bounded and
$\|\phi\|_{cb}=\sup_n\|\phi_n\|$ is finite.

\section{ Dilations of Operator Valued Measures}

Possibly the first well-known dilation result for operator valued measures is due to Naimark.

\begin{theorem}[Naimark's Dilation Theorem]\label{th:23}
Let $E$ be a regular, positive, $B(\cH)$-valued measure on $\Omega$.
Then there exist a Hilbert space $\cK$, a bounded linear operator
$V:\cH\to\cK$, and a regular, self-adjoint, spectral,
$B(\cK)$-valued measure $F$ on $\Omega$, such that
\[E(B)=V^*F(B)V.\]
\end{theorem}

From Naimark's dilation
Theorem, we know that every regular positive operator-valued measure (OVM for abbreviation) can be
dilated to a self-adjoint, spectral operator-valued measure on a
larger Hilbert space. But not all of the operator-valued measures
can have a Hilbert dilation space.  Such an example was constructed in \cite{HLLL1} in which we constructed an operator-valued
measure induced by the framing that does not have a Hilbert
dilation space. The construction is based on an example of Osaka \cite{Os} of a normal
non-completely bounded map from $\ell^{\infty}(\Bbb{N})$ into $B(H)$. In fact operator valued measures that admit Hilbert space dilations are the ones that are closely related to completely
bounded measures and maps.

Now let $\Omega$ be a compact Hausdorff space, let $E$ be a bounded,
regular, operator-valued measure on $\Omega,$ and let
$\phi:C(\Omega)\to B(\cH)$ be the bounded, linear map associated
with $E$ by integration. So for any
$f\in C(\Omega),$
\[\langle \phi(f)x,y \rangle=\int_\Omega f\, d\, \mu_{x,y},\]
where
\[ \mu_{x,y}(B)=\langle E(B)x,y \rangle \]
The OVM $E$ is called completely bounded when $\phi$ is completely
bounded. Using Wittstock's decomposition theorem, $E$ is completely
bounded if and only if it can be expressed as a linear combination
of positive operator-valued measures.

Let $\{x_i\}_{i\in\J}$ be a non-zero frame for a separable Hilbert
space $\cH.$ Let $\Sigma$ be the $\sigma$-algebra of all subsets of
$\J$, and
$$E: \Sigma\to B(\cH), \quad
E(B)=\sum_{i\in B} x_i\otimes x_i$$
Since  $E$ is a regular, positive $B(\cH)$-valued measure, by
Naimark's dilation Theorem \ref{th:23}, there exists a Hilbert space
$\cK$, a bounded linear operator $V:\cH\to\cK$, and a regular,
self-adjoint, spectral, $B(\cK)$-valued measure $F$ on $\J$, such
that
\[E(B)=V^*F(B)V.\]
This Hilbert space $\cK$ can be $\ell_2$, and the atoms
$x_{i}\otimes x_{i}$ of the measure dilate to rank-1 projections
$e_{i}\otimes e_{i}$, where $\{e_{i}\}$ is the standard orthonormal
basis for $\ell^{2}$.  That is $\cK$ can be the same as the dilation
space in Proposition \ref{prop:orthog} (ii).

In the case that  $\{x_i,y_i\}_{i\in\J}$ is a non-zero framing
for a separable Hilbert space $\cH, $ and
$E(B)=\sum_{i\in B} x_i\otimes y_i$ for all $B\in \Sigma$,
$E$ is a $B(\cH)$-valued measure. In \cite{HLLL1} we showed that this $E$
also has a dilation space $Z$. But this dilation space is not
necessarily a Hilbert space, in general, it is a Banach space and
consistent with Proposition \ref{pr:33}. The dilation is essentially
constructed using Proposition \ref{pr:33} (ii), where the dilation
of the atoms $x_{i}\otimes y_{i}$ corresponds to the projection
$u_{i}\otimes u_{i}^{*}$ and $\{u_{i}\}$ is an unconditional basis
for the dilation space $Z$.

Framings are the natural generalization of discrete frame theory
(more specifically, dual-frame pairs) to non-Hilbertian settings.
Even if the underlying space is a Hilbert space,  the dilation space
for framing induced operator valued measures can fail to be
Hilbertian. This theory was originally developed by Casazza, Han and Larson in \cite{CHL} as an attempt to introduce {\it frame theory
    with dilations} into a Banach space context.  The initial motivation of this investigation was to completely understand
    the dilation theory of framings.  In the context of Hilbert spaces, we realized that the dilation theory for discrete
    framings from \cite{CHL} induces a dilation theory for discrete operator valued measures that may fail to be completely bounded.

These examples inspired us
to consider Banach space dilation theory for arbitrary operator valued measures.

\begin{definition}\label{de:35}
Let $X$ and $Y$ be Banach spaces, and let $(\Omega,\Sigma)$ be a
measurable space. A \emph{$B(X,Y)$-valued measure on $\Omega$} is a
map $E:\Sigma\to B(X,Y)$ that is countably additive in the weak
operator topology; that is, if $\{B_i\}$ is a disjoint countable
collection of members of $\Sigma$ with union $B$, then
\begin{equation*}
y^{*}(E(B)x)=\sum_i y^{*} (E(B_i)x)
\end{equation*}
for all $x\in X$ and $y^*\in Y^*.$
\end{definition}

We will use the symbol $(\Omega, \Sigma, E)$ if the range space is
clear from context, or $(\Omega, \Sigma, E, B(X, Y))$, to denote
this operator-valued measure system.

The Orlicz-Pettis theorem  states that weak unconditional
convergence and norm unconditional convergence of a series are the
same in every Banach space (c.f. \cite{DJT}). Thus  we have  that
$\sum_i E(B_i)x$ weakly unconditionally converges to $E(B)x$ if
and only if $\sum_i E(B_i)x$ strongly unconditionally converges to
$E(B)x.$ So Definition \ref{de:35} is equivalent to saying that
$E$ is strongly countably additive, that is, if $\{B_i\}$ is a
disjoint countable collection of members of $\Sigma$ with union
$B$, then
\begin{equation*}
E(B)x=\sum_i E(B_i)x,\quad  \ \forall x\in X.
\end{equation*}

\begin{definition}\label{de:37}
Let $E$ be a $B(X,Y)$-valued measure on $(\Omega,\Sigma).$ Then the \textit{norm of $E$} is defined by
\[\|E\|=\sup_{B\in\Sigma}\|E(B)\|.\]
We call $E$ \textit{normalized} if $\|E\|=1.$
\end{definition}

A $B(X,Y)$-valued measure $E$ is always bounded, i.e.
\begin{equation}\label{eq:31}
\sup_{B\in\Sigma}\|E(B)\|<+\infty.
\end{equation}
Indeed, for all $x\in X$ and $y^*\in Y^*,$
$\mu_{x,y^*}(B):=y^*(E(B)x) $ is a complex measure on
$(\Omega,\Sigma).$ From complex measure theory (c.f. \cite{Rud}),
we know that $\mu_{x,y^*}$ is bounded, i.e.
\[\sup_{B\in\Sigma}|y^*(E(B)x)|<+\infty.
\]
By the Uniform Boundedness Principle, we get (\ref{eq:31}).

Similar to the Hilbert space operator valued measures, we introduce the following definitions.

\begin{definition}\label{de:7}
A $B(X)$-valued measure $E$ on $(\Omega,\Sigma)$ is called:
\begin{enumerate}
\item[(i)] an \textit{operator-valued probability measur}e if
$E(\Omega)=\mathrm{I}_X,$

\item[(ii)]a \textit{projection-valued measure} if $E(B)$ is a projection
on $X$ for all $B\in\Sigma$,

\item[(iii)] a \textit{spectral operator-valued measure} if for all $A,
B\in \Sigma, E(A\cap B)=E(A)\cdot E(B)$ (we will also use the term
\textit{idempotent-valued measure} to mean a spectral-valued measure.)
\end{enumerate}
\end{definition}

For general operator valued measures we established the following dilation theorem \cite{HLLL1}.

\begin{theorem} Let $E:\Sigma \to B(X,Y)$ be an operator-valued measure. Then there exist a Banach space $Z$,
 bounded linear operators $S:Z\to Y$ and $T:X\to Z$, and a projection-valued probability measure
 $F:\Sigma\to B(Z)$ such that $$E(B)=SF(B)T$$ for all $B\in\Sigma$.
 \end{theorem}

\noindent  We will call $(F, Z, S, T)$ in the above theorem a {\it
Banach space dilation system}, and a {\it Hilbert dilation system}
if $Z$ can be taken as a Hilbert space. This theorem  generalizes
Naimark's (Neumark's) Dilation Theorem for positive operator valued
measures.  But even in the case that the underlying space is a
Hilbert space the dilation space cannot always be taken to be a
Hilbert space.  Thus elements of the theory of Banach spaces are
essential in this work.

A key idea is the introduction of the
elementary dilation space  and  and the minimal dilation norm.

Let $X,Y$ be Banach spaces and $(\Omega,\Sigma,E,B(X,Y))$  an
operator-valued measure system. For any $B\in \Sigma$ and $x\in X,$
define
\[E_{ B, x }:\Sigma\to Y,\quad E_{ B, x }(A)=E(B\cap A)x,
\quad \forall A\in \Sigma.\]
 Then it is easy to see that
$E_{ B, x }$ is a vector-valued measure on $(\Omega,\Sigma)$ of
$Y$ and $E_{ B, x }\in \mathfrak{M}^Y_\Sigma.$

Let $M_E=\mbox{span}\{E_{ B, x }: x\in X, B\in \Sigma\}.$  We introduce some linear mappings on the spaces $X,$ $Y$ and $M_E.$

For any $\{C_i\}_{i=1}^N\subset
\mathbb{C},$ $\{B_i\}_{i=1}^N\subset \Sigma$ and
$\{x_i\}_{i=1}^N\subset X$, the mappings
$$S:M_E\to Y,\quad
S\Big(\sum_{i=1}^NC_i E_{{B_i}, {x_i}}\Big)=\sum_{i=1}^NC_i
E(B_i)x_i$$
$$T: X\to M_E,\quad T(x)=E_{\Omega,x}$$
and
$$F(B): M_E\to M_E,\quad F(B)\Big(\sum_{i=1}^NC_i E_{ {B_i}, {x_i}}\Big)
=\sum_{i=1}^NC_i E_{{B\cap B_i}, {x_i}}, \quad \forall B\in\Sigma$$
are well-defined and linear.

\begin{definition}\label{de:418}
Let $M_E$ be the space induced by $(\Omega,\Sigma,E,B(X,Y)).$ Let
$\|\cdot\|$ be a norm on $M_E$. Denote this normed space by
$M_{E,\|\cdot\|}$ and its completion $\widetilde{M}_{E,\|\cdot\|}.$
The norm on $\widetilde{M}_{E,\|\cdot\|}$, with $\|\cdot\| :=
\|\cdot\|_{\cD}$ given by a norming function $\cD$ as discussed
above, is called a \textit{dilation norm of $E$} if the following
conditions are satisfied:
\begin{enumerate}
\item[(i)] The mapping $S_\cD: \widetilde{M}_{E,\cD}\to Y$ defined on $M_E$ by
\begin{eqnarray*}
S_\cD\Big(\sum_{i=1}^NC_i E_{ {B_i}, {x_i}}\Big)=\sum_{i=1}^NC_i E(B_i)x_i
\end{eqnarray*}
is bounded. \item[(ii)]
 The mapping $T_\cD: X\to \widetilde{M}_{E,\cD}$ defined by
\begin{eqnarray*}
T_\cD(x)=E_{\Omega,x}
\end{eqnarray*}
is bounded. \item[(iii)]
 The mapping $F_\cD: \Sigma\to B(\widetilde{M}_{E,\cD})$ defined
 by
\begin{eqnarray*}
F_\cD(B)\Big(\sum_{i=1}^NC_i E_{ {B_i}, {x_i}}\Big)
=\sum_{i=1}^NC_i E_{{B\cap B_i},{x_i}}
\end{eqnarray*}
is an operator-valued measure,
\end{enumerate}
where $\{C_i\}_{i=1}^N\subset \mathbb{C},$ $\{x_i\}_{i=1}^N\subset
X$ and $\{B_i\}_{i=1}^N\subset \Sigma$.
\end{definition}

We call the Banach space
$\widetilde{M}_{E,\cD}$ \textit{the elementary dilation space of $E$} and
$$(\Omega,\Sigma,F_\cD,B(\widetilde{M}_{E,\cD}),S_\cD,T_\cD)$$ the elementary dilation operator-valued measure system.
 The \textit{minimal dilation norm}
$\|\cdot\|_\alpha$ on $M_E$ is defined by
\[\Big\|\sum^{N}_{i=1}C_iE_{ {B_i}, {x_i}}\Big\|_\alpha
=\sup_{B\in\Sigma}\Big\|\sum_{i=1}^N C_iE(B\cap B_i)x_i\Big\|_Y\]
for all $\sum^{N}_{i=1}C_iE_{ {B_i}, {x_i}}\in M_E$. Using this we
show that every OVM has a projection valued dilation to an
elementary dilation space, and moreover,  $||\cdot ||_{\alpha}$ is a
minimal norm on the elementary dilation space.

A corresponding dilation projection-valued measure system
$(\Omega,\Sigma,F,B(Z),S,T)$ is said to be \textit{injectiv}e if
$\sum F(B_i)T(x_i)=0$ whenever $\sum E_{B_i,x_i}=0$.

It is useful to note that all the elementary dilation spaces are
Banach spaces of functions.

\begin{theorem} Let $E:\Sigma\to B(X,Y)$ be an operator-valued measure and $(F,Z,S,T)$ be an injective Banach
 space dilation system.  Then we have the following:

(i) There exist an elementary Banach space dilation system $(F_\cD,
\widetilde{M}_{E,\cD},S_\cD,T_\cD)$ of $E$ and a linear isometric
embedding
\[U:\widetilde{M}_{E,\cD}\to Z\]
such that
\[S_\mathcal {D}=SU,\ F(\Omega)T=UT_\cD,\ UF_\cD(B)=F(B)U,
\quad \forall B\in \Sigma.\]

(ii) The norm $\|\cdot\|_\alpha$ is indeed a dilation norm.
Moreover,  If $\cD$ is a dilation norm of $E,$ then there exists a
constant $C_\cD$ such that for any $\sum_{i=1}^NC_iE_{ {B_i},
{x_i}}\in M_{E,\cD},$
\[
\sup_{B\in\Sigma}\Big\|\sum_{i=1}^N C_iE(B\cap B_i)x_i\Big\|_Y
\leq C_\cD \Big\|\sum_{i=1}^NC_iE_{ {B_i}, {x_i}}\Big\|_\cD,
\]
where $N>0,$ $\{C_i\}_{i=1}^{N}\subset\mathbb{C},$
$\{x_i\}_{i=1}^{N}\subset X$ and $\{B_i\}_{i=1}^{N}\subset \Sigma. $
Consequently
\[ \|f\|_\alpha\le C_\cD\|f\|_\cD, \qquad \forall f\in M_{E}.\]
\end{theorem}
\bigskip




\begin{definition}
Let $E:\Sigma\to B(X,Y)$ be an operator-valued measure and $(F,Z,S,T)$ be a Banach
space dilation system. Then $(F,Z,S,T)$ is called {\it linearly minimal} if
$Z$ is the closed linear span of $F(\Sigma)TX$, where $F(\Sigma)TX=\{ F(B)(T x) : B\in\Sigma, x\in X \}$.
\end{definition}

A projection valued measure can have a nontrivial linearly minimal
dilation to another projection valued measure.  The following simple
example illustrates this.  It is also an example of a dilation
projection-valued measure system which is not injective and for
which the conclusion of Theorem 3.7 is not true. This shows that if
we drop the ``injectivity" in the hypothesis of Theorem 3.7 the
conclusion need not be true. However a simple modification of the
conclusion will be true (see Remark 3.10).

%
%

\begin{example} Let $(\Omega, \Sigma, \mu)$ be a probability space and let
$\nu$ be a finite measure that dominates $\mu$.
Let $X := L^2(\Omega, \mu)$ and let $Y := L^2(\Omega, \nu)$. Let $
\alpha$ be a bounded linear functional on $X$ that takes  $1$ at the
function $\eta  =  1$. Let $\Omega = \Omega_{0}^{c} \cup \Omega_0$
be the Hahn
 decomposition, where
$ \Omega_0$ is a measurable subset of $\Omega$ which is a null set for $\mu$
and which
 supports  the singular part of $\nu$ with respect to $\mu$. Regard $L^2(\Omega, \nu)$  as the
direct sum of $L^2(\Omega, \mu)$ and $L^2(\Omega_0, \nu)$.
 Embed $X$ into $Y$ by $T(f) = f \oplus \alpha(f)\chi_{\Omega_0
}$, where $\chi_{\Omega_0 }$ is the constant function $1$ in
$L^2(\Omega_0, \nu)$.  Since $\alpha$ is a linear functional $T$ is
a linear map.  In particular it maps the constant function $1$ in $X
:= L^2(\Omega, \mu)$ to the constant function $1$ in $Y :=
L^2(\Omega, \nu)$.
Define a projection valued measure $\phi : \Sigma \rightarrow B(X)$
by setting $\phi(B) = M_{\chi_{B}}$, the projection operator of
multiplication by the characteristic function of $B$. Do the same
construction to define a projection valued measure $\Phi: \Sigma
\rightarrow B(Y$).  Since $TX$ contains the constant function $1$ in
$L^2(\Omega, \nu)$, the closed linear span of $\Phi(\Sigma)TX$ is
$Y$.

Let $S$ denote the mapping of $Y := L^2(\Omega, \nu)$ onto $X :=
L^2(\Omega, \mu)$ determined by the function mapping $f \rightarrow
f|_{{\Omega_0}^c}$. Then  $S$ has kernel $ L^2(\Omega_0, \nu)$.

Then $\Phi$ is a dilation of $\phi$ for the dilation maps $T$ and
$S$, and the dilation is linearly minimal because the closed linear
span of $\Phi(\Sigma)TX$ is $Y$.  The dilation is clearly
non-injective, and the conclusion of Theorem 3.7 fails for it.
\end{example}


\begin{remark}
We have the following natural generalization of Theorem 3.7: Let all
terms be as in the hypotheses of Theorem 3.7 except do not assume
that the Banach space dilation system $(F,Z,S,T)$ is injective.
First, obtain a reduction if necessary by restricting the range
space of $F$ so that the closure of the range of $F$ times the range
of $T$ is all of $Z$. This makes the dilation linearly minimal.
Example 3.9 shows that this reduction to linearly minimal is not
alone  sufficient to generalize Theorem 3.7.  Obtain a second
reduction by replacing $Z$ with its quotient by the kernel of $S$.
Then the hypotheses of Theorem 3.7 are satisfied, so we can obtain a
generalization of Theorem 3.7 by removing the injectivity
requirement in the hypothesis and inserting the restriction
reduction followed by the quotient reduction in the statement of the
conclusion. In Example 3.9 the restriction reduction is unnecessary
because the dilation is already linearly minimal, and the quotient
reduction makes the reduced dilation equivalent to $\phi$.

The point of this is that the {\it minimal elementary norm dilation}
of this section is really a {\it geometrically minimal dilation} in
the sense that any dilation, after a simple restriction reduction
and a quotient reduction if necessary, is isometrically isomorphic
to an {\it elementary} dilation norm dilation.  And the class of
elementary dilation norm spaces are related in the sense that there
is a minimal dilation norm and a maximal dilation norm, and all
dilation norms lie between the minimal and the maximal norm on the
elementary function space, and the actual dilation space is the
completion of the elementary function space in one of the dilation
norms. So in this sense the minimal norm elementary dilation of an
operator valued measure is subordinate to all other dilations of the
OVM.
\end{remark}


 While in general an operator-valued probability
    measure does not admit a Hilbert space dilation, the dilation
    theory can be strengthened in the case that it does admit a Hilbert space dilation:

\begin{theorem} Let $E:\Sigma\to B(\cH)$ be an operator-valued probability measure.
 If $E$ has a Hilbert dilation system $(\widetilde{E},\widetilde{H},S,T)$, then there exists a
corresponding Hilbert dilation system $(F,\cK, V^*, V)$ such that
$V:\cH\to\cK$ is an isometric embedding.

\end{theorem}

\noindent   This theorem turns out to have some interesting
applications to framing induced operator valued measure dilation. In
particular,  it led to a complete characterization of framings whose
induced operator valued measures are completely bounded. We include
here a few sample examples with the following theorem:

\begin{theorem} \label{thm-rescale} Let $(x_i,y_i)_{i\in\N}$ be a non-zero framing for a Hilbert space
$\cH$, and  $E$ be the operator-valued probability measure induced
by $(x_i,y_i)_{i\in\N}$. Then we have the following:

(i) $E$ has a Hilbert dilation space $\cK$ if and only if there
exist $\alpha_i,\beta_i \in \C, i\in\N$ with
$\alpha_i\bar{\beta_i}=1$ such that $\{\alpha_ix_i\}_{i\in\N}$ and
$\{\beta_iy_i\}_{i\in\N}$ both are the frames for the Hilbert space
$\cH.$

(ii)  $E$ is a completely bounded map if and only if
$\{x_i,y_i\}_{i\in\N}$ can be re-scaled to dual frames.

(iii)   If $\ \inf\|x_i\|\cdot\|y_i\|>0,$ then we can find
$\alpha_i,\beta_i \in \C, i\in\N$ with $\alpha_i\bar{\beta_i}=1$
such that $\{\alpha_ix_i\}_{i\in\N}$ and $\{\beta_iy_i\}_{i\in\N}$
both are frames for the Hilbert space $\cH$. Hence the
operator-valued measure induced by $\{x_i,y_i\}_{i\in\N}$ has a
Hilbertian dilation.

\end{theorem}

\noindent For the existence of non-rescalable (to dual frame pairs)
framings, we obtained the following:

\begin{theorem} \label{thm-example} There exists a framing for a Hilbert space such that its induced operator-valued measure is not completely bounded,
and consequently it can not be re-scaled to obtain a framing that
admits a Hilbert space dilation.
\end{theorem}

The second part of this theorem follows from the first part of Theorem  \ref{thm-rescale} (ii).

\begin{remark} For the existence of such an example, the motivating example
of framing constructed by Casazza, Han and Larson (Example 3.9 in
\cite{CHL}) can not be dilated  to an unconditional basis for a Hilbert
space, although it can be dilated to an unconditional basis for a
Banach space. We originally conjectured that this is an example that
fails to induce a completely bounded operator valued measure.
However, it turns out that this framing can be re-scaled to a
framing that admits a Hilbert space dilation , and consequently
disproves our conjecture. Our construction of the new example in
Theorem \ref{thm-example} uses a non-completely bounded map to construct a
non-completely bounded OVM which yields the required framing. This
delimiting example shows that the dilation theory for framings
developed in \cite{CHL} gives a true generalization of Naimark's
Dilation Theorem for the discrete case. This is the example that led
us to consider general (non-necessarily-discrete) operator valued
measures, and to the results of Chapter 2 that lead to the dilation
theory for general (not necessarily completely bounded) OVM's that
completely generalizes Naimark's Dilation theorem in a Banach space
setting, and which is new even for Hilbert spaces.
\end{remark}

Part (iii) of Theorem  \ref{thm-rescale} provides us a sufficient condition under
which a framing induced operator-valued measure has a Hilbert space
dilation. This can be applied to framings that have nice structures.
For example, the following is an unexpected result for unitary
system induced framings, where a unitary system is a countable
collection of unitary operators. This clearly applies to wavelet and
Gabor systems.

\begin{corollary} Let $\mathcal{U}_1$ and $\mathcal{U}_2$ be unitary systems on a
separable Hilbert space $\cH.$ If there exist $x,y\in\cH$ such that
$\{\mathcal{U}_1x,\mathcal{U}_2y\}$ is a framing of $\cH,$ then
$\{\mathcal{U}_1x\}$ and $\{\mathcal{U}_2y\}$ both are frames for
$\cH.$
\end{corollary}

There exist examples of sequences $\{x_{n}\}$ and $\{y_{n}\}$
in a Hilbert space $\mathcal H$ with the following properties:

(i) $x = \sum_{n} \langle x, x_{n}\rangle y_{n}$  hold for all $x$
in a dense subset of $\mathcal H$, and the convergence is
unconditionally.

(ii) $\inf ||x_{n}||\cdot ||y_{n}|| >0$.

(iii) $\{x_{n},  y_{n}\}$ is not a framing.

\begin{example} Let $\mathcal{H} = L^{2}[0,
1]$, and $g(t) = t^{1/4}$, $f(t) = 1/g(t)$. Define $x_{n}(t) =
e^{2\pi i nt}f(t)$ and $y_{n}(t) = e^{2\pi in t}g(t)$. Then it is
easy to verify $(i)$ and $(ii)$. For $(iii)$, we consider the
convergence of the series
$$
\sum_{n\in \Bbb{Z}}\langle f, x_{n}\rangle y_{n}.
$$
Note that $||\langle f, x_{n}\rangle y_{n}||^{2} = |\langle f,
x_{n}\rangle|^{2} \cdot ||g||^{2}$ and $\{\langle f,
x_{n}\rangle\}$ is not in $\ell^{2}$ (since $f^{2} \notin L^{2}[0,
1]$). Thus
$\sum_{n}\langle f, x_{n}\rangle y_{n}$ can not be convergent
unconditionally. Therefore $\{x_{n},  y_{n}\}$ is not a framing.
\end{example}

\section{Dilations of Bounded Linear Maps}

Inspired by the techniques used to build the dilation theory for
general operator valued measures we consider establishing a dilation
theory for general linear maps. Historically the dilation theory has
been extensively investigated in the context of positive, or
completely bounded maps on C*-algebras, with Stinespring's dilation
theorem as possibly one of the most notable results in this
direction (c.f. \cite{Arv, Pa} and the references therein).


\begin{theorem}\label{Stines-dilation} [Stinespring's dilation theorem]\label{th:22}
Let $\mathcal {A}$ be a unital $C^*$-algebra, and let $\phi:\mathcal
{A}\to B(\cH)$ be a completely positive map. Then there exists a
Hilbert space $\cK,$ a unital $*-$homomorphism $\pi: \mathcal {A}\to
B(\cK),$ and a bounded operator $V:\cH\to\cK$ with
$\|\phi(1)\|=\|V\|^2$ such that
\[\phi(a)=V^*\pi(a)V.\]
\end{theorem}

 The following is also well known for commutative  $C^*$-algebras:

\begin{theorem}[cf. Theorem 3.11, \cite{Pa}]\label{th:21}
Let $\cB$ be a $C^*$-algebra, and let $\phi:C(\Omega)\to\cB$ be
positive. Then $\phi$ is completely positive.
\end{theorem}

This result together with Theorem \ref{th:22} implies that
Stinespring's dilation theorem holds for positive maps when
$\mathcal{A}$ is a unital commutative $C^*$-algebra.

A proof of Naimark's dilation theorem by using Stinespring's
dilation theorem can be sketched as follows: Let $\phi:{\mathcal
A}\to B(\cH)$ be the natural extension of $E$ to the $C^*$-algebra
${\mathcal A}$ generated by all the characteristic functions of
measurable subsets of $\Omega$.  Then $\phi$ is positive, and hence
is completely positive by Theorem \ref{th:21}. Apply Stinespring's
dilation theorem to obtain a $*-$homomorphism $\pi: {\mathcal A} \to
B(\cK),$ and a bounded, linear operator $V:\cH\to\cK$ such that
$\phi(f)=V^*\pi(f)V$ for all $f$ in ${\mathcal A}$.  Let $F$ be the
$B(\cK)-$valued measure corresponding to $\pi.$ Then it can be
verified that $F$ has the desired properties.

Completely positive maps are completely bounded. In the other
direction we have Wittstock's decomposition theorem \cite{Pa}:

\begin{proposition}
Let $\mathcal {A}$ be a unital $C^*$-algebra, and let $\phi:\mathcal
{A}\to B(\cH)$ be a completely bounded map. Then $\phi$ is a linear
combination of two completely positive maps.
\end{proposition}

The following  is a generalization of Stinespring's representation
theorem.

\begin{theorem}\label{de:25}
Let $\mathcal {A}$ be a  unital $C^*$-algebra, and let
$\phi:\mathcal {A}\to B(\cH)$ be a completely bounded map. Then
there exists a Hilbert space $\cK,$ a $*-$homomorphism $\pi:
\mathcal {A}\to B(\cK),$ and bounded operators $V_i:\cH\to\cK,
i=1,2,$ with $\|\phi\|_{cb}=\|V_1\|\cdot\|V_2\|$ such that
\[\phi(a)=V_1^*\pi(a)V_2\]
for all $a\in\mathcal {A}.$ Moreover, if $\|\phi\|_{cb}=1,$ then
$V_1$ and $V_2$ may be taken to be isometries.
\end{theorem}

Now let $\Omega$ be a compact Hausdorff space, let $E$ be a bounded,
regular, operator-valued measure on $\Omega,$ and let
$\phi:C(\Omega)\to B(\cH)$ be the bounded, linear map associated
with $E$ by integration as described in section 1.4.1. So for any
$f\in C(\Omega),$
\[\langle \phi(f)x,y \rangle=\int_\Omega f\, d\, \mu_{x,y},\]
where
\[ \mu_{x,y}(B)=\langle E(B)x,y \rangle. \]

The OVM $E$ is called completely bounded when $\phi$ is completely
bounded. Using Wittstock's decomposition theorem, $E$ is completely
bounded if and only if it can be expressed as a linear combination
of positive operator-valued measures.

One of the important applications  of  our main dilation theorem is
the  dilation for not necessarily cb-maps with appropriate
continuity properties from a commutative von Neumann algebra into
$B(\cH)$. While the ultraweak topology on $B(\cH)$ for a Hilbert space
$\cH$ is well-understood,  we define the ultraweak topology on
$B(X)$ for a Banach space $X$ through tensor products: Let
$X\otimes Y$ be the tensor product of the Banach space $X$ and
$Y.$ The projective norm on $X\otimes Y$ is defined by:
\[\|u\|_{\wedge}=\inf\Big\{\sum_{i=1}^n\|x_i\|\|y_i\|:u=\sum_{i=1}^n x_i\otimes y_i\Big\}.\]
We will use $X\otimes_{\wedge} Y$ to denote the tensor product
$X\otimes Y$ endowed with the projective norm
$\|\cdot\|_{\wedge}.$   Its completion will be denoted by
$X\widehat{\otimes} Y.$ From \cite{R} Section 2.2, for any Banach
spaces $X$ and $Y,$ we have the identification:
\[(X\widehat{\otimes} Y)^*=B(X,Y^*).\]
Thus $B(X,X^{**})=(X\widehat{\otimes}  X^*)^*.$ Viewing $X\subseteq
X^{**},$ we define the {\it ultraweak topology} on $B(X)$ to be the
weak* topology induced by the predual $X\widehat{\otimes} X^*.$  We
will also use the term {\it normal} to denote an ultraweakly
continuous linear map.

\begin{theorem} If $\mathcal{A}$ is a purely
atomic abelian von Neumann algebra acting on a separable Hilbert
space, then for every ultraweakly continuous linear map
$\phi:\mathcal{A}\to B(\cH)$, there exists a Banach space $Z,$ an
ultraweakly continuous unital homomorphism $\pi:\mathcal{A}\to
B(Z)$, and bounded linear operators $T:\cH\to Z$ and $S:Z\to\cH$
such that $$\phi(a)=S\pi(a)T$$ for all $a\in\mathcal{A}$.
\end{theorem}

\noindent The proof of this theorem uses some special properties of the
minimal dilation system for the $\phi$ induced operator valued
measure on the space $(\Bbb{N}, 2^{\Bbb{N}})$. Motivated by some
ideas used in the proof of the above theorem,  we then obtained a {\it universal
dilation theorem} for all bounded linear mappings between Banach algebras:

\begin{theorem} \label{thm-universal} Let $\mathcal{A}$ be a Banach
algebra, and let $\phi:\mathcal{A}\to B(X)$ be a bounded linear
operator, where $X$ is a Banach space. Then there exists a Banach
space $Z,$ a bounded linear unital homomorphism $\pi:\mathcal{A}\to
B(Z)$, and bounded linear operators $T:X\to Z$ and $S:Z\to X$
such that
\[\phi(a)=S\pi(a)T\]
for all $a\in\mathcal{A}.$
\end{theorem}

Since this theorem is so general we would expect that there is a also purely algebraic dilation theorem
for any linear transformations. This indeed is the case.

\begin{proposition}\label{th:t50}
If $A$ is unital algebra, $V$  a vector space, and
$\phi:A\rightarrow L(V)$ a linear map, then there exists a vector
space $W$, a unital homomorphism $\pi:A\rightarrow L(V)$, and linear
maps $T:V\rightarrow W$, $S:W\rightarrow V$, such that
$$\phi(\cdot)=S\pi(\cdot) T.$$
\end{proposition}

This result  maybe well-known. However we provide a short proof for interested readers.

\begin{proof}
For $a\in A, x\in V$, define $\alpha_{a,x}\in L(A,V)$ by
$$\alpha_{a,x}(\cdot):=\phi(\cdot a) x.$$ Let $W:=\mathrm{span}\{\alpha_{a,x}:a\in A,x\in V\}\subset L(A,V).$
Define $\pi:A\rightarrow L(W)$ by
$\pi(a)(\alpha_{b,x}):=\alpha_{ab,x}.$ It is easy to see that
$\pi$ is a unital homomorphism. For $x\in V$ define $T:V
\rightarrow L(A,V)$ by $T_x:=\alpha_{I,x}=\phi(\cdot
I)x=\phi(\cdot)x$. Define $S:W\rightarrow W$ by setting
$S(\alpha_{a,x}):=\phi(a)x$ and extending linearly to $W$. If
$a\in A, x\in V$ are arbitrary, we have
$S\pi(a)Tx=S\pi(a)\alpha_{I,x}=S \alpha_{a,x}=\phi(a)x.$ Hence
$\phi=S\pi T.$
\end{proof}

We note that the above proposition has been generalized by the
second author and F. Szafraniec \cite{LS} to the case where $A$ is a
unital semigroup.

Theorem \ref{thm-universal}  is a true generalization of our commutative
theorem in an important special case, and generalizes some of our
results for maps of commutative von Neumann algebras to the case
where the von Neumann algebra is non-commutative.

For the case when
$\mathcal{A}$ is a von Neumann algebra acting on a separable Hilbert
space and $\phi$ is ultraweakly continuous (i.e., normal) we
conjecture that  the dilation space $Z$ can be taken to be separable
and the dilation homomorphism $\pi$ is also ultraweakly continuous.
While we are not able to confirm this conjecture we have the
following result. Here, SOT is the abbreviation of strong operator topology.

\begin{theorem} Let $K,H$ be Hilbert spaces,
$A\subset B(K)$ be a von Neumann algebra, and $\phi:A\rightarrow
B(H)$ be a bounded linear operator which is ultraweakly-\textsc{SOT}
continuous on the unit ball $B_A$ of $A$. Then there exists a Banach
space $Z$, a bounded linear homeomorphism $\pi:A\rightarrow B(Z)$
which is \textsc{SOT}-\textsc{SOT} continuous on $B_A$, and bounded
linear operator $T:H\rightarrow Z$ and $S:Z\rightarrow H$ such that
$$\phi(a)=S \pi(a)T$$ for all $a\in A.$ If in addition that $K,H$ are separable, then
the Banach space $Z$ can be taken to be separable.
\end{theorem}

These results are apparently new for mappings of von Neumann
algebras. They generalize special cases of Stinespring's Dilation
Theorem. The standard discrete Hilbert space frame theory is
identified with the special case of our theory in which the domain
algebra is abelian and purely atomic, the map is completely bounded,
and the OVM is purely atomic and completely bounded with rank-1
atoms.

The universal  dilation result has connections with Kadison's
similarity problem for bounded homomorphisms between von Neumann
algebras (see the Remark 4.14). For example,  if $\mathcal{A}$
belongs to one of the following classes:  nuclear; $\mathcal{A} =
B(H)$; $\mathcal{A}$ has no tracial states;  $\mathcal{A}$ is
commutative; $II_{1}$-factor with Murry and von Neumann's property
$\Gamma$, then any non completely bounded map $\phi:\mathcal{A}\to
B(H)$  can never have a Hilbertian dilation (i.e.  the dilation
space $Z$ can never be a Hilbert space)  since otherwise
$\pi:\mathcal{A}\to B(Z)$ would be similar to a *-homomorphism and
hence completely bounded and so would be $\phi$. On the other hand,
if there exists a  von Neumann algebra $\mathcal{A}$  and a non
completely bounded map $\phi$ from $\mathcal{A}$ to $B(H)$ that has
a Hilbert space dilation: $\pi:\mathcal{A}\to B(Z)$ (i.e., where $Z$
is a Hilbert space), then $\pi$ will be a counterexample to the
Kadison's similarity problem since in this case $\pi$ is a
homomorphisim that is not completely bounded and consequently can
not be similar to a *-homomorphisim.

\section{Some Remarks and Problems}

\begin{remark} It is well known that there is a  theory establishing a connection
between general bounded linear mappings from the $C^*$-algebra
$C(X)$ of continuous functions on a compact Hausdorf space $X$ into
$B(H)$ and operator valued measures on the sigma algebra of Borel
    subsets of $X$ (c.f. \cite{Pa}).  If $A$ is an abelian $C^*$-algebra then $A$ can be identified with $C(X)$ for a topological space $X$
     and can also be identified with $C(\beta X)$ where  $\beta X$ is the Stone-Cech compactification of X.  Then the
     support $\sigma$-algebra for the OVM is the sigma algebra of Borel subsets of $\beta X$  which is enormous.  However in our generalized
     (commutative) framing theory $\mathcal{A}$ will always be an abelian von Neumann algebra presented up front
     as $L^\infty(\Omega, \Sigma, \mu) $, with $\Omega$ a topological space and $\Sigma$ its algebra of Borel sets,
     and the maps on $A$ into $B(H)$ are  normal.  In particular, to model the discrete  frame and framing theory $\Omega$
     is a countable index set with the discrete topology (most often $\Bbb{N}$), so $\Sigma$ is its power set, and $\mu$ is
     counting measure.   So in this setting it is more natural to work directly with this presentation in developing
     dilation theory rather than passing to $\beta \Omega$, and we took this approach in our investigation.
        \end{remark}

        \begin{remark}

     We feel that the connection we make with established discrete
frame and framing theory is transparent,
    and then the OVM dilation theory for the continuous case becomes a natural but nontrivial generalization of the theory
    for the discrete case that was inspired by framings. After doing this we attempted to apply our techniques to the case
    where the domain algebra
    for a map is non-commutative.    However, additional hypotheses are needed
    if dilations of maps are to have strong
    continuity and structural properties.  For a map between C*-algebras it is well-known that there is a Hilbert space dilation
    if the map is
    completely bounded. (If the domain algebra is commutative this statement is an iff.) Even if a map is not cb it has a Banach
    space dilation.
     We are interested in the continuity and structural properties a dilation can have.  In  the
     discrete abelian case,
 the dilation of a normal map can be taken to be normal and the dilation space can be taken to be separable, and with suitable hypotheses this type of result
 can be generalized to the noncommutative setting.
\end{remark}

The following is a list of problems we think may be important for
the general dilation theory of operator valued measures and bounded
linear maps.

It was proven in \cite{HLLL1} that if  $\{i\}$ is an atom in
$\Sigma$ and $E$ is an operator valued frame on $\Sigma$, then the
minimal dilation $F_{\alpha}$ has the property that the rank of
$F_\alpha(\{i\})$ is equal to the rank of  $E(\{i\})$. This leads to
the following problem.
\bigskip

\noindent{\bf Problem 1.} Is it always true that with an appropriate
notion of rank function for an operator valued measure, that $r
(F_{\alpha} (B)) = r( E (B))$ for every $B\in \Sigma$? What about if
a ``rank" definition is defined by: $r(B) = \sup \{\mathrm{rank}
E(A): A\subset B, A\in\Sigma\}$?

\bigskip

Let $(\Omega,\Sigma, \mu)$ be a  probability space and let $\phi:
L^\infty(\mu)\rightarrow B(\mathcal{H})$ be ultraweakly continuous.
Then it naturally induces an operator valued probablity measure
$$
E(B) = \phi(\chi_{B}), \ \ \forall \,B\in \Sigma.
$$

\bigskip

\noindent{\bf Problem 2.} Let $E: (\Omega,\Sigma)\rightarrow
B(\mathcal{H})$ be an operator valued measure. Is there an
ultraweakly continuous map $\phi:L^\infty(\mu)\rightarrow
B(\mathcal{H})$ that induces $E$ on $(\Omega,\Sigma)?$ If the answer
is negative, then determine necessary and sufficient conditions for
$E$ to be induced by an ultraweakly continuous map?
\bigskip

 As with Stinespring's dilation theorem, if $A$ and $\cH$ in Theorem \ref{Stines-dilation}  are both separable then
 the dilated Banach
 space $Z$ is also
separable. However, the Banach algebras  we are interested in include von Neumann algebras and these
 are generally not separable, and the linear maps $\phi:A\rightarrow B(H)$ are often normal. So we pose the
following two problems.
\bigskip

\noindent{\bf Problem 3.} Let $K,H$ be separable Hilbert spaces, let
$A\subset B(K)$ be a von Neumann algebra, and let $\phi:A\rightarrow
B(H)$ be a bounded linear map. When is there a {\it separable}
Banach space $Z,$ a bounded linear unital homomorphism
$\pi:\mathcal{A}\to B(Z)$, and bounded linear operators $T:\cH\to Z$
and $S:Z\to\cH$ such that
\[\phi(a)=S\pi(a)T\]
for all $a\in\mathcal{A}$ ?

\bigskip

\noindent{\bf Problem 4.} Let  $A\subset B(K)$ be a von Neumann
algebra, and $\phi:A\rightarrow B(H)$ be a normal linear map. When
can we dilate $\phi$ to a {\it normal} linear unital homomorphism
$\pi:\mathcal{A}\to B(Z)$ for some (reflexive) Banach space $Z$?

\bigskip

Finally, concerning the Hilbert space dilations and  Kadison's
Similarity Problem, we are interested in the following questions:
\bigskip

\noindent{\bf Problem 5.} Let  $A\subset B(K)$ be a von Neumann
algebra, and let $\phi:A\rightarrow B(H)$ be a bounded linear map.
We know that $\phi$ has a Hilbert space dilation if it is completely
bounded. Is there a non-completely bounded map that admits a Hilbert
space dilation? In particular, if $\phi(A) = A^{t}$ for any $A \in
\oplus_{n=1}^{\infty} M_{n\times n}(\C)$, then $T$ is bounded but
not completely bounded. What can we say about the dilation of
$\phi$? Does it admit a Hilbert space dilation?
\bigskip

Yes. An affirmative answer would yield a negative answer to the similarity problem.
\bigskip

\noindent{\bf Problem 6.}  Let  $A\subset B(K)$ be a von Neumann
algebra, and let $\phi:A\rightarrow B(H)$ be a bounded linear map.
``Characterize"  those maps that admit Hilbert space dilations, and
``Characterize" those maps that admit reflexive Banach space
dilations.

\end{document}